\magnification=\magstep1


\def\item{\vskip1.3pt\hang\textindent}


\tolerance=300
\pretolerance=200
\hfuzz=1pt
\vfuzz=1pt

\hoffset 0cm            
\hsize=5.8 true in
\vsize=9.5 true in

\def\rightheadline{\hfil\smc\lastname\hfil\tenbf\folio}
\def\leftheadline{\tenbf\folio\hfil\smc\lastname\hfil}
\headline={\ifodd\pageno\rightheadline\else\leftheadline\fi}
\newdimen\dimenone
\def\checkleftspace#1#2#3#4#5{
 \dimenone=\pagetotal
 \advance\dimenone by -\pageshrink   
 \ifdim\dimenone>\pagegoal          
   \else\dimenone=\pagetotal
        \advance\dimenone by \pagestretch
        \ifdim\dimenone<\pagegoal
          \dimenone=\pagetotal
          \advance\dimenone by#1         
          \setbox0=\vbox{#2\parskip=0pt                
                       \hyphenpenalty=10000
                       \rightskip=0pt plus 5em
                       \noindent#3 \vskip#4}    
        \advance\dimenone by\ht0
        \advance\dimenone by 3\baselineskip
        \ifdim\dimenone>\pagegoal\vfill\eject\fi
          \else\eject\fi\fi}

\parindent=35pt
\mathsurround=1pt
\parskip=1pt plus .25pt minus .25pt
\normallineskiplimit=.99pt

\mathchardef\emptyset="001F 

\def\Res{\mathop{\rm Res}\limits}
\def\Int{\mathop{\rm int}\nolimits}
%



\def\1{{\bf1}}\def\0{{\bf0}}

\def\({\bigl(}  \def\){\bigr)}
\def\<{\mathopen{\langle}}\def\>{\mathclose{\rangle}}

\def\Z{{\mathchoice{{\hbox{$\rm Z\hskip 0.26em\llap{\rm Z}$}}}%
{{\hbox{$\rm Z\hskip 0.26em\llap{\rm Z}$}}}%
{{\hbox{$\scriptstyle\rm Z\hskip 0.31em\llap{$\scriptstyle\rm Z$}$}}}{{%
\hbox{$\scriptscriptstyle\rm Z$\hskip0.18em\llap{$\scriptscriptstyle\rm Z$}}}}}}

\def\N{{\mathchoice{\hbox{$\rm I\hskip-0.14em N$}}%
{\hbox{$\rm I\hskip-0.14em N$}}%
{\hbox{$\scriptstyle\rm I\hskip-0.14em N$}}%
{\hbox{$\scriptscriptstyle\rm I\hskip-0.10em N$}}}}

\def\R{{\mathchoice{\hbox{$\rm I\hskip-0.14em R$}}%
{\hbox{$\rm I\hskip-0.14em R$}}%
{\hbox{$\scriptstyle\rm I\hskip-0.14em R$}}%
{\hbox{$\scriptscriptstyle\rm I\hskip-0.10em R$}}}}

\def\C{{\mathchoice%
{\hbox{$\rm C\hskip-0.47em\hbox{%
\vrule height 0.58em width 0.06em depth-0.035em}$}\;}%
{\hbox{$\rm C\hskip-0.47em\hbox{%
\vrule height 0.58em width 0.06em depth-0.035em}$}\;}%
{\hbox{$\scriptstyle\rm C\hskip-0.46em\hbox{%
$\scriptstyle\vrule height 0.365em width 0.05em depth-0.025em$}$}\;}
{\hbox{$\scriptscriptstyle\rm C\hskip-0.41em\hbox{
$\scriptscriptstyle\vrule height 0.285em width 0.04em depth-0.018em$}$}\;}}}

\def\.{{\cdot}}
\def\|{\Vert}
\def\ssk{\smallskip}
\def\msk{\medskip}
\def\bsk{\bigskip}
\def\giantskip{\vskip2\bigskipamount}

\def\giantbreak{\par \ifdim\lastskip<2\bigskipamount \removelastskip
         \penalty-400 \giantskip\fi}

\def\nin{\noindent}
\def\cen{\centerline}
\def\pagebreak{\vskip 0pt plus 0.0001fil\break}
\def\linebreak{\break}

\def\epsilon{\varepsilon}

\font\ninerm=cmr9
\font\eightrm=cmr8
\font\sixrm=cmr6

\font\eightbf=cmbx8
\font\sixbf=cmbx6

\font\eighti=cmmi8
\font\sixi=cmmi6
\font\ninesy=cmsy9
\font\eightsy=cmsy8
\font\sixsy=cmsy6

\font\eightit=cmti8


\font\eightsl=cmsl8

\font\eighttt=cmtt8
\font\bfone=cmbx10 scaled\magstep1 
\font\smc=cmcsc10

\font\small=cmcsc8

\def\no #1. {\bigbreak\vskip-\parskip\noindent\bf #1. \quad\rm}

\def\Proposition #1. {\checkleftspace{0pt}{\bf}{Theorem}{0pt}{}
\bigbreak\vskip-\parskip\noindent{\bf Proposition #1.}
\quad\it}

\def\Theorem #1. {\checkleftspace{0pt}{\bf}{Theorem}{0pt}{}
\bigbreak\vskip-\parskip\noindent{\bf  Theorem #1.}
\quad\it}
\def\Corollary #1. {\checkleftspace{0pt}{\bf}{Theorem}{0pt}{}
\bigbreak\vskip-\parskip\nin{\bf Corollary #1.}
\quad\it}
\def\Lemma #1. {\checkleftspace{0pt}{\bf}{Theorem}{0pt}{}
\bigbreak\vskip-\parskip\noindent{\bf  Lemma #1.}\quad\it}

\def\Definition #1. {\checkleftspace{0pt}{\bf}{Theorem}{0pt}{}
\rm\bigbreak\vskip-\parskip\noindent{\bf Definition #1.}
\quad}

\def\Remark #1. {\checkleftspace{0pt}{\bf}{Theorem}{0pt}{}
\rm\bigbreak\vskip-\parskip\noindent{\bf Remark #1.}\quad}

\def\Exercise #1. {\checkleftspace{0pt}{\bf}{Theorem}{0pt}{}
\rm\bigbreak\vskip-\parskip\noindent{\bf Exercise #1.}
\quad}

\def\Example #1. {\checkleftspace{0pt}{\bf}{Theorem}{0pt}{}
\rm\bigbreak\vskip-\parskip\noindent{\bf Example #1.}\quad}
\def\Examples #1. {\checkleftspace{0pt}{\bf}{Theorem}{0pt}
\rm\bigbreak\vskip-\parskip\noindent{\bf Examples #1.}\quad}

\newcount\problemnumb \problemnumb=0
\def\Problem{\global\advance\problemnumb by 1\bigbreak\vskip-\parskip\noindent
{\bf Problem \the\problemnumb.}\quad\rm }

\def\Proof#1.{\rm\par\ifdim\lastskip<\bigskipamount\removelastskip\fi\smallskip
            \noindent {\bf Proof.}\quad}

\nopagenumbers

\def\author{}
\def\lastname{}
\def\thanks#1{\footnote*{\eightrm#1}}
\def\title{}

\def\nonumbers{\def\leftheadline{\hfil} \def\rightheadline{\hfil}}

\def\lastname{}
\def\h{{\textstyle{1\over2}}} 
 
\def\he{{1\over2}} 
\def\si{\sigma}

\def\ep{\epsilon}

\def\text{\textstyle} 
\def\disp{\displaystyle} 
\def\d{{\,\rm d}}

\def\and{{\rm and }}

\expandafter\edef\csname amssym.def\endcsname{%
       \catcode`\noexpand\@=\the\catcode`\@\space}
\catcode`\@=11
\def\undefine#1{\let#1\undefined}
\def\newsymbol#1#2#3#4#5{\let\next@\relax
 \ifnum#2=\@ne\let\next@\msafam@\else
 \ifnum#2=\tw@\let\next@\msbfam@\fi\fi
 \mathchardef#1="#3\next@#4#5}
\def\mathhexbox@#1#2#3{\relax
 \ifmmode\mathpalette{}{\m@th\mathchar"#1#2#3}%
 \else\leavevmode\hbox{$\m@th\mathchar"#1#2#3$}\fi}
\def\hexnumber@#1{\ifcase#1 0\or 1\or 2\or 3\or 4\or 5\or 6\or 7\or 8\or
 9\or A\or B\or C\or D\or E\or F\fi}

\font\tenmsb=msbm10
\font\sevenmsb=msbm7
\font\fivemsb=msbm5
\newfam\msbfam
\textfont\msbfam=\tenmsb
\scriptfont\msbfam=\sevenmsb
\scriptscriptfont\msbfam=\fivemsb
\edef\msbfam@{\hexnumber@\msbfam}
\def\Bbb#1{{\fam\msbfam\relax#1}}

\newsymbol\Bbbk 207C
\def\widehat#1{\setbox\z@\hbox{$\m@th#1$}%
 \ifdim\wd\z@>\tw@ em\mathaccent"0\msbfam@5B{#1}%
 \else\mathaccent"0362{#1}\fi}
\def\widetilde#1{\setbox\z@\hbox{$\m@th#1$}%
 \ifdim\wd\z@>\tw@ em\mathaccent"0\msbfam@5D{#1}%
 \else\mathaccent"0365{#1}\fi}
\font\teneufm=eufm10
\font\seveneufm=eufm7
\font\fiveeufm=eufm5
\newfam\eufmfam
\textfont\eufmfam=\teneufm
\scriptfont\eufmfam=\seveneufm
\scriptscriptfont\eufmfam=\fiveeufm

\catcode`@=11 

\expandafter\edef\csname amssym.def\endcsname{%
       \catcode`\noexpand\@=\the\catcode`\@\space}
\font\eightmsb=msbm8
\font\sixmsb=msbm6
\font\fivemsb=msbm5
\font\eighteufm=eufm8
\font\sixeufm=eufm6
\font\fiveeufm=eufm5
\newskip\ttglue
\def\eightpoint{\def\rm{\fam0\eightrm}%
  \textfont0=\eightrm \scriptfont0=\sixrm \scriptscriptfont0=\fiverm
  \textfont1=\eighti \scriptfont1=\sixi \scriptscriptfont1=\fivei
  \textfont2=\eightsy \scriptfont2=\sixsy \scriptscriptfont2=\fivesy
  \textfont3=\tenex \scriptfont3=\tenex \scriptscriptfont3=\tenex
\textfont\eufmfam=\eighteufm
\scriptfont\eufmfam=\sixeufm
\scriptscriptfont\eufmfam=\fiveeufm
\textfont\msbfam=\eightmsb
\scriptfont\msbfam=\sixmsb
\scriptscriptfont\msbfam=\fivemsb
  \def\it{\fam\itfam\eightit}%
  \textfont\itfam=\eightit
  \def\sl{\fam\slfam\eightsl}%
  \textfont\slfam=\eightsl
  \def\bf{\fam\bffam\eightbf}%
  \textfont\bffam=\eightbf \scriptfont\bffam=\sixbf
   \scriptscriptfont\bffam=\fivebf
  \def\tt{\fam\ttfam\eighttt}%
  \textfont\ttfam=\eighttt
  \tt \ttglue=.5em plus.25em minus.15em
  \normalbaselineskip=9pt
  \def\MF{{\manual opqr}\-{\manual stuq}}%
  \let\big=\eightbig
  \setbox\strutbox=\hbox{\vrule height7pt depth2pt width\z@}%
  \normalbaselines\rm}
\def\eightbig#1{{\hbox{$\textfont0=\ninerm\textfont2=\ninesy
  \left#1\vbox to6.5pt{}\right.\n@space$}}}


\csname amssym.def\endcsname


\def\la{\lambda} 
\def\al{\alpha} 
\def\be{\beta}

\def\({\left(} 
\def\){\right)} 
 
\def\eq{\eqalign} 
 
\def\gcd{{\rm gcd}} 
\def\O#1{O\(#1\)} 
\def\abs#1{\left| #1 \right|}

\def\klein{\eightpoint \def\smc{\small} \baselineskip=9pt}   

\def\fn#1#2{{\parindent=0.7true cm
\footnote{$^{(#1)}$}{{\klein  #2}}}}

\font\boldmas=msbm10                  
\def\Bbb#1{\hbox{\boldmas #1}}        
\def\Z{{\Bbb Z}}                        
\def\N{{\Bbb N}}                        

\def\R{{\Bbb R}}

\def\C{{\Bbb C}}


\font\eightrm=cmr8                                                    
\long\def\fussnote#1#2{{\baselineskip=9pt                            
\setbox\strutbox=\hbox{\vrule height 7pt depth 2pt width 0pt}%
\eightrm                                                         
\footnote{#1}{#2}}}                                              
\font\boldmasi=msbm10 scaled 700      
\def\Bbbi#1{\hbox{\boldmasi #1}}      
\font\boldmas=msbm10                  
\def\Bbb#1{\hbox{\boldmas #1}}        
\def\Zi{{\Bbbi Z}}                      
\def\Pi{{\Bbbi P}}                      
\def\Ri{{\Bbbi R}}



\def\dint #1 {
\quad  \setbox0=\hbox{$\disp\int\!\!\!\int$}
  \setbox1=\hbox{$\!\!\!_{#1}$}
  \vtop{\hsize=\wd1\centerline{\copy0}\copy1} \quad}

\def\drint #1 {
\qquad  \setbox0=\hbox{$\disp\int\!\!\!\int\!\!\!\int$}
  \setbox1=\hbox{$\!\!\!_{#1}$}
  \vtop{\hsize=\wd1\centerline{\copy0}\copy1}\qquad}

\def\frac#1#2{{#1\over #2}}

\def\date{\the\day.~\the\month.~\the\year}

\def\klein{\eightpoint \def\smc{\small} }

\def\frac#1#2{{#1\over#2}} 
\def\Int{\int\limits}

\nonumbers
\hsize=16true cm     \vsize=24true cm

\parindent=0cm 

\def\m{(m,n)} 
\def\zz{\Z_*^2}
\def\zzi{\Zi_*^2}
\def\z{z_0^*} 
\def\b{\be_0^*} 
\def\T{{\cal T}_Q}

\vbox{\vskip0.5true cm} 
\cen{\bfone Primitive lattice points inside an ellipse} \bsk 
\cen{\bf Werner Georg NOWAK} \bsk


\vbox{\vskip0.7true cm} 

{\bf 1.~Introduction.}\quad Let $Q=Q\m=am^2+bmn+cn^2$ be a positive definite 
binary quadratic form, where $a,b,c$ are 
{\it arbitrary} real numbers with $a>0,\ D:=4ac-b^2>0$. For a large parameter 
$x$, we consider the lattice point quantities 
$$  \eq{A(x) = &\ \# \{\m\in\Z^2:\ Q\m \le x\}\,,    \cr 
B(x) = &\ \# \{\m\in\Z^2:\ Q\m \le x,\ \gcd\m=1\ \}\,,  \cr  } \eqno(1.1) $$ 
which count the number of {\it all}, resp., of all {\it primitive} lattice 
points in the ellipse disc $Q \le x$. It is well known that 
$$  A(x) = {2\pi\over\sqrt{D}} x + P(x)\,,\qquad 
B(x) = {12\over\pi\sqrt{D}} x + R(x)\,,  \eqno(1.2) $$ 
where $P(x), R(x)$ are error terms on which a lot of research has been done. 
(For an enlightening presentation of this theory, see the monograph of 
Kr\"atzel [10].) As far as $P(x)$ is concerned, the sharpest 
published\fn{1}{Actually, M.~Huxley has meanwhile improved further his upper 
bound, essentially replacing the exponent ${23\over73}=0.315068\dots$ 
by ${131\over416}=0.314903\dots$. 
The author is indebted to Professor Huxley for sending him a copy of his 
unpublished manuscript.} results read 
$$  P(x) \ll x^{23/73} \(\log x\)^{315/146}\,,  \eqno(1.3) $$ 
$$  \liminf_{x\to\infty}\({P(x)\over x^{1/4} (\log x)^{1/4} }\) < 0\,, 
\eqno(1.4) $$ 
and 
$$   \int_0^T (P(t^2))^2 \d t \sim C_Q\, T^2\,. \eqno(1.5) $$ 
They are due to M.~Huxley [6], [7], the author [14], P.~Bleher [1] 
and the author [15]\fn{2}{In this latter reference, 
actually a short interval version of this 
asymptotics is established. We omit the discussion of a possible error term 
in (1.5) which, for the case of a general ellipse, is by no means simple.}. 
All these estimates have been proved for general 
convex planar domains with smooth boundary of nonvanishing curvature. \msk 
The question for analogous results about $R(x)$ remains much more enigmatic. 
To see why, we recall that the generating Dirichlet series corresponding to 
$P(x)$, resp., $A(x)$, is the Epstein zeta-function 
$$  \zeta_Q(s) = \sum_{\m\in\zzi} Q\m^{-s} \qquad(\Re(s)>1)\,, \eqno(1.6)  $$ 
where $\zz:=\Z^2\setminus\{(0,0)\}$. It possesses an analytic continuation to 
the whole complex plane, with the exception of a simple pole at $s=1$, 
and satisfies a functional equation 
$$  \zeta_Q(s) = \({2\pi\over\sqrt{D}}\)^{2s-1} \, 
{\Gamma(1-s)\over\Gamma(s)}\,\zeta_Q(1-s)\,.  \eqno(1.7)  $$  
(See Potter [17], or, for a multivariate version, Kr\"atzel's monograph [11], 
p.~202.) By Vinogradov's Lemma, the generating function of $B(x)$ reads, 
for $\Re(s)>1$, 
$$ \sum_{\m\in\zzi\atop\gcd\m=1} Q\m^{-s} = \sum_{k=1}^\infty \mu(k)
\sum_{\m\in\zzi} Q(km,kn)^{-s} ={\zeta_Q(s)\over\zeta(2s)}\,. \eqno(1.8) $$ 
By Perron's formula, for every value of $x>0$ which is not attained by $Q\m$, 
$\m\in\zz$, 
$$  B(x) = {1\over2\pi i} \Int_{2-i\infty}^{2+i\infty} 
{\zeta_Q(s)\over\zeta(2s)}\,{x^s\over s} \d s\,.  $$ 
Shifting the line of integration to the left, we are confronted with the lack 
of information about the zeros of the Riemann 
zeta-function\fn{3}{For an enlightening presentation of its theory the reader 
is referred to the monograph of \hbox{A.~Ivi\'c [9].}}: These might come close 
to $\Re(s)=1$, hence an estimate $R(x)\ll x^\theta$ cannot be proved for any 
$\theta<\h$, at the present state of art. The best known upper bound is 
$$  R(x) = \O{x^{1/2} \exp\(-C(\log x)^{3/5} (\log\log x)^{-1/5}\)}\,.  $$ 
Several authors have investigated this problem under the assumption of the 
Riemann Hypothesis. After previous work by Huxley \& Nowak [8] and by 
W.~M\"uller [13], the sharpest {\it conditional } results of this kind 
are due to W.~Zhai [24] and read $R(x)\ll x^{221/608+\ep}$ for a rational form 
$Q$, and $R(x)\ll x^{33349/84040+\ep}$ in general. 
(Note that ${221\over608}=0.3634\dots$, ${33349\over84040}=0.3968\dots$ .) 
See also Zhai \& Cao [23] and Wu [22]. 
\msk There is little hope to establish estimates for $R(x)$ which are directly 
analogous to (1.4) and (1.5). \ssk 
Nevertheless, in the present paper we shall prove a result which says that at 
least the lower bound part of (1.5) holds true for $R(x)$ 
also.\fn{4}{Ironically, our analysis actually will yield this result not 
{\it although} there is the cumbersome denominator $\zeta(2s)$ in (1.8), but 
{\it because} it is there.} Trivially, 
this implies a pointwise $\Omega$-result for $R(x)$, which is comparable to, 
though slightly weaker than, formula (1.4). \bsk 

\vbox{{\bf Theorem.}\quad{\it The error term $R(x)$ defined in $(1.1), (1.2)$ 
satisfies 
$$  {1\over Y} \Int_1^Y \abs{R(x)} \d x \gg Y^{1/4}\,,  \eqno(1.9)  $$ 
as $Y\to\infty$, the $\gg$-constant depending on the form $Q$.}}

\bsk\msk  

{\bf 2.~A zero-density bound for Epstein zeta-functions.}\fn{5}{The result 
stated suffices for our purpose and will be believed at first glance by the 
expert. However, it is difficult to find it explicitly in the literature. 
Further, it cannot be improved substantially: As Davenport \& Heilbronn [3], 
[4], and M.~Voronin [20] showed, if $Q$ is an integral form of 
class number exceeding 1, then 
$N_Q(1,T)\gg T$ and also $N_Q(\al,T)-N_Q(1,T)\gg T$ for $\h<\al<1$.} \bsk 

{\bf Lemma.}\quad{\it For any positive definite binary quadratic form $Q$, \ 
$\si\in\R$ and $T\in\R^+$, denote by $ N_Q^*(\si,T)$ the number of zeros 
(counted with multiplicity) of $\zeta_Q(s)$ with $\Re(s)=\si,\ 
\abs{\Im(s)}\le T$, and put $N_Q(\si,T)=\sum\limits_{\si'\ge\si}N_Q^*(\si',T)$. 
Then, as $T\to\infty$,  
$$  N_Q^*({\text{1\over4}},T) = N_Q^*({\text{3\over4}},T) \le 
N_Q({\text{3\over4}},T) = o(T\log T)\,.  $$ } \bsk 

{\bf Proof.}\quad First of all, 
$N_Q^*({\text{1\over4}},T) = N_Q^*({\text{3\over4}},T)$ is clear by the 
functional equation (1.7). To establish the $o$-assertion, one can follow the 
classical example of Titchmarsh's monograph [19], section 9.15. We rewrite 
(1.6), 
for $\Re(s)>1$, as 
$$ \zeta_Q(s) = \sum_{k=1}^\infty r_k \la_k^{-s} = 
\la_1^{-s}\(r_1 + U(s)\)\,, $$ 
where $r_k\in\N^*$ and $\(\la_k\)$ is a strictly increasing sequence of 
positive reals. Since $U(\si+it)\to0$ as $\si\to\infty$, uniformly in 
$t$, there exists some $\sigma^*>1$ (depending on $Q$) such that 
$\abs{U(\si+it)}\le\h r_1$ for $\si\ge\si^*$ and all $t$. As a consequence, 
$$   \abs{\zeta_Q(\si^* +it)}\ge \h r_1 \la_1^{-\si^*}  \eqno(2.1)  $$ 
for all $t$, and $\zeta_Q(s)\ne0$ for $\Re(s)\ge\si^*$. Let further 
$\T := \{t\in\R:\ \cos(t \log \la_1)\ge{\text{3\over4}}\}$, then 
$$ \abs{\Re(\zeta_Q(\si^*+it))} \ge {\text{1\over4}} r_1 \la_1^{-\si^*}  
\eqno(2.2)  $$ for all $t\in\T$.  \ssk  
We use a variant of 
formula (9.9.1) in [19] ("Littlewood's Lemma"): 
If $\al>0$ and $T>0,\ T\in\T$\fn{6}{Obviously, for any given $T_0\in\Ri^+$, 
there exists some $T\in\T$ with $T_0\le T \ll T_0$.} 
are such that there are no zeros of 
$\zeta_Q(s)$ on $\Re(s)=\al$ and on $\abs{\Im(s)}=T$, then 
$$  \Int_{{\cal R}} \log\zeta_Q(s) \d s = 
-2\pi i \Int_\al^{\si^*} N_Q(\si,T)\d\si + O(1)\,,  \eqno(2.3)  $$ 
where ${\cal R}$ is the rectangle $(\al\pm iT), (\si^*\pm i T)$, 
and the logarithm is defined (almost everywhere) by 
$$\log\zeta_Q(\si+it)=\log\zeta_Q(\si^*) 
+ \Int_{{\cal C}} {\zeta_Q'(s)\over\zeta_Q(s)} \d s$$ where 
$\log\zeta_Q(\si^*)\in\R$ and ${\cal C}$ consists 
of the two straight line segments from $\si^*$ to $\si^*+it$ and further to 
$\si+it$. Moreover, let ${\rm arg}\,\zeta_Q(s) := \Im(\log\zeta_Q(s))$. \par 
Taking the imaginary part of (2.3), we get 
$$ \eq{2\pi \Int_\al^{\si^*} N_Q(\si,T)\d\si = & \Int_{-T}^T 
\log\abs{\zeta_Q(\al+it)}\d t - 
\Int_{-T}^T \log\abs{\zeta_Q(\si^*+it)}\d t + \cr & + \Int_\al^{\si^*} 
{\rm arg}\,\zeta_Q(\si+iT)\d\si - \Int_\al^{\si^*} 
{\rm arg}\,\zeta_Q(\si-iT)\d\si + O(1)\,. \cr  }   $$ 
By (2.1), the second integral on the right hand side is $O(T)$.  
We mimick the argument in section 9.4 of [19] to show that (at least) 
$$ {\rm arg}\,\zeta_Q(\si\pm iT)=O(T)\eqno(2.4)$$ 
uniformly in $\al\le\si\le\si^*$.                              
This will readily yield 
$$ 2\pi \Int_\al^{\si^*} N_Q(\si,T)\d\si = 
\Int_{-T}^T \log\abs{\zeta_Q(\al+it)}\d t + O(T)\,,  \eqno(2.5)  $$ 
for any fixed $\al>0$ and $T\to\infty$. To prove (2.4), we note first that 
$\zeta_Q'(s)\over\zeta_Q(s)$ is bounded on $\Re(s)=\si^*$, hence 
$\arg\zeta_Q(\si^*\pm iT) = O(T)$. The variation of 
$\arg\zeta_Q(\si \pm iT)$ on $\al\le\si\le\si^*$ is $\ll 1+q$, $q$ the number 
of zeros of $\Re(\zeta_Q(\si \pm iT))$ on this line segment. Further, $q\le 
n(\si^*-\al)$, if $n(r)$ denotes the number of zeros (counted with 
multiplicity) of the function $G(s):=\h(\zeta_Q(s\pm iT) + \zeta_Q(s\mp iT))$ 
in the disc $\abs{s-\si^*}\le r$. Now 
$$ \Int_0^{\si^*-\he\al} {n(r)\over r} \d r \ge 
\Int_{\si^*-\al}^{\si^*-\he\al} {n(r)\over r} \d r \gg n(\si^*-\al)\,, $$ 
and, by Jensen's theorem, 
$$  \Int_0^{\si^*-\he\al} {n(r)\over r} \d r = 
{1\over2\pi}\Int_0^{2\pi} \log\abs{G(\si^*
+(\si^*-\h\al)e^{i\theta})}\d\theta - \log\abs{G(\si^*)} \ll \log T\,, $$ 
since $\abs{G(\si^*)}\gg1$ because of $T\in\T$ and (2.2). This establishes 
(2.4) and thus (2.5). \ssk 
According to W.~M\"uller 
[13]\fn{7}{In fact, M\"uller proves this bound more generally for the Hlawka 
zeta-function of a convex planar domain with smooth boundary of nonvanishing 
curvature. Similar results can be found in Huxley \& Nowak [8] and in 
W.~Zhai [24].}, Prop.~2, at least for every $\al\ge{\text{2\over3}}$, 
$$ \Int_0^T \abs{\zeta_Q\(\al+it\)}^2\d t \ll T^{1+\ep}  $$ 
for any $\ep>0$. Hence, 
by Jensen's inequality (e.g., [5], p.~1132) and the reflection principle, 
for suitable $\al\in]{\text{2\over3}},{\text{3\over4}}[$,  
$$  \Int_{-T}^T \log\abs{\zeta_Q(\al+it)}\d t 
\le T \log\({1\over T} \Int_0^T \abs{\zeta_Q(\al+it)}^2 \d t\) \ll 
\ep T \log T\,.    $$ 
Thus, by (2.5), for $\si_0=\h(\al+{\text{3\over4}})$, 
$$ N_Q(\si_0,T) \le {1\over\si_0-\al} \Int_{\al}^{\si_0} N_Q(\si,T)\d\si \ll 
\ep T\log T\,.  $$ 
Since $\ep>0$ is arbitrary, this establishes the lemma. \bsk\msk 

\vbox{{\bf 3.~Proof of the Theorem.}\quad 
Following an idea due to Pintz [16], we consider the Mellin transform, for 
$\Re(s)>1$, 
$$  \eq{H(s):=& \Int_1^\infty R(x) x^{-s-1} \d x = 
\Int_1^\infty \(\sum_{Q\m\le x\atop\gcd\m=1} 1 - {12\over\pi\sqrt{D}} x\) 
x^{-s-1} \d x = \cr = & \sum_{\m\in\zzi\atop\gcd\m=1} \Int_{Q\m}^\infty 
x^{-1-s} \d x - {12\over\pi\sqrt{D}} \Int_1^\infty x^{-s} \d x =\cr 
= & {\zeta_Q(s)\over s\,\zeta(2s)} - {12\over\pi\sqrt{D}}\,{1\over s-1} =: 
{E(s)\over s(s-1)\zeta(2s)(2s-1)}\,.  \cr   } \eqno(3.1) $$ }

Obviously $H(s)$ possesses a meromorphic continuation to all of $\C$, with  
$E(s)$ an entire function. Now choose $z_0={\text{1\over4}}+i\be_0$ such that 
$2z_0$ is a zero of the Riemann zeta-function 
and $\zeta_Q(z_0)\ne0$. (The existence follows from the above 
lemma and a celebrated result of Selberg [18], refined further by Levinson [12] 
and Conrey [2].) The function 
$$ g(s) := {s(s-1)\zeta(2s)(2s-1)\over(s-z_0)(s+2)^7} \eqno(3.2) $$ 
is regular in $\Re(s)>-2$, and so is 
$$ g(s)H(s) = {E(s)\over(s-z_0)(s+2)^7}\,, \eqno(3.3) $$ 
apart from a simple pole at $s=z_0$, since 
$E(z_0)=(z_0-1)(2z_0-1)\zeta_Q(z_0)\ne0$. 
By the functional equation (1.7), 
$\zeta_Q(-1+it) \ll \abs{t}^3$, 
and similarly $\zeta(-2+2it)\asymp\abs{t}^{5/2}$, 
as $\abs{t}\to\infty$, hence the integrals 
$\Int_{\be-i\infty}^{\be+i\infty}\abs{g(s)} \d s$ and 
$\Int_{\be-i\infty}^{\be+i\infty}\abs{g(s)\,H(s)} \d s$ 
converge for $\be\in\{-1,2\}$. For $\eta>0$, we define a weight function 
$$   w(\eta) := \Int_{2-i\infty}^{2+i\infty} g(s) \eta^{s+1} \d s\,, 
\eqno(3.4) $$ 
which satisfies 
$$  w(\eta) = \cases{O(1) & for $\eta\ge1$, 
\cr 0 & for $0<\eta<1$.} \eqno(3.5) $$ 
(To see this, one can shift the line of integration to 
$\Int_{-1-i\infty}^{-1+i\infty}$ in the first case and to 
$\Int_{C-i\infty}^{C+i\infty}$, with $C\to\infty$, in the second case.) 
Thus, for $Y>0$, 

\vbox{$$ V(Y) := {1\over Y}\Int_1^\infty R(x) w\({Y\over x}\) \d x = 
{1\over Y}\Int_1^\infty R(x) 
\(\,\Int_{2-i\infty}^{2+i\infty} g(s) \({Y\over x}\)^{s+1} \d s\) \d x = $$ 
$$ = \Int_{2-i\infty}^{2+i\infty} g(s) Y^s \( \Int_1^\infty R(x) x^{-s-1} 
\d x\) \d s = \Int_{2-i\infty}^{2+i\infty} g(s) H(s) Y^s \d s\,. \eqno(3.6) 
$$}
 
Shifting the line of integration to $\Re(s)=-1$, we get, for $Y$ large, 
$$  V(Y) = 2\pi i\,\Res_{s=z_0}\(g(s)H(s)Y^s\) + 
\Int_{-1-i\infty}^{-1+i\infty} g(s) H(s) Y^s \d s = 2\pi i\,\al_0 Y^{z_0} 
+ O(Y^{-1})\,, \eqno(3.7) $$ 
where 
$$  \al_0 = {E(z_0)\over(z_0+2)^7} 
= {(z_0-1)(2z_0-1)\zeta_Q(z_0)\over(z_0+2)^7}\,.  \eqno(3.8)  $$  
From this it is evident that, as $Y\to\infty$,  
$$   \abs{V(Y)} \gg \abs{Y^{z_0}} = Y^{1/4} \eqno(3.9) $$ 

\vbox{and, on the other hand, in view of (3.5), 
$$  \abs{V(Y)} = \abs{{1\over Y}\Int_1^Y R(x) w\({Y\over x}\)\d x} 
\ll {1\over Y}\Int_1^Y \abs{R(x)}\d x\,, \eqno(3.10) $$ 
which completes the proof of our theorem.} \bsk\msk

{\bf 4.~How to get an estimate with an explicit constant.}\quad The above 
argument was clearly non-effective, as far as the $\gg$-constant in (1.9) is 
concerned: In particular, our lemma only guarantees the existence of a 
Riemann-zeta zero $2z_0$ for which $\zeta_Q(z_0)\ne0$, but gives no possibility 
to estimate it. \par 
In this final section, we shall therefore show how to obtain a lower 
bound\fn{8}{However, we shall not invest too much effort to make this bound 
as large as possible.} for 
$$ K_0 := \liminf_{Y\to\infty}\(Y^{-5/4} \Int_1^Y \abs{R(x)} \d x\)\,,$$ 
for any specific given form $Q\m$. Our first step is to show that 
$$  \abs{w(\eta)}  \le 0.33\,, \eqno(4.1)  $$ 
for all $\eta>0$ and any $Q$. In fact, by (3.4) and (3.2), 
$$  \eq{\abs{w(\eta)} = & 
\abs{\Int_{-1-i\infty}^{-1+i\infty} g(s) \eta^{s+1} \d s} 
\le \Int_{-\infty}^{\infty} \abs{g(-1+it)} \d t  \le \cr 
\le & \Int_{-\infty}^{\infty} \abs{(-1+it)(-2+it)(-
3+2it)\over(1+it)^7 (-{5\over4}+i(t-\be_0))}\abs{\zeta(-2+2it)} \d t\,, \cr }$$ 
if we recall that $z_0={1\over4}+i\be_0$. We further 
use the functional equation (e.g., [19], f.~(2.1.9), (2.1.10)) 
$$   \zeta(-2+2it) = \pi^{-5/2 +2it} \, 
{\Gamma({3\over2}-it)\over\Gamma(-1+it)} \zeta(3-2it)\,,  $$ 
along with well-known identities for the $\Gamma$-function (in particular 
f.~8.332 in [5]) which imply 
$$ \abs{\Gamma({3\over2}-it)\over\Gamma(-1+it)} \le 
\abs{(\h+it)(-1+it)}\sqrt{\abs{t}}\,. $$ 
Thus 
$$  \abs{w(\eta)} \le {\zeta(3) \over \pi^{5/2}}\,
\Int_{-\infty}^\infty \abs{(-1+it)^2(-2+it)(-3+2it)(\h+it)\over
(1+it)^7 (-{5\over4}+i(t-\be_0))}\,\sqrt{\abs{t}}\, \d t\ \le $$ 
$$ \le {\zeta(3) \over \pi^{5/2}}\,\(2 \Int_0^\infty 
{(4+t^2)(9+4t^2)({1\over4}+t^2) t \over(1+t^2)^5}\,\d t\, 
\Int_{-\infty}^\infty {\d t\over {25\over16}+(t-\be_0)^2}\,\)^{1/2}\,,  $$ 
by Cauchy's inequality. The integrals are evaluated to $143\over32$ (with a 
little help from {\it Mathematica } [21], e.g.) and $4\pi\over5$, which 
readily gives (4.1).  
By (3.10) and (3.7), it follows that 
$$  K_0 \ge 6\pi\,\abs{\al_0}\,, \eqno(4.2) $$ 
thus it remains to estimate $\abs{\al_0}$ (see (3.8)), in particular 
$\abs{\zeta_Q(z_0)}$, for any fixed form $Q$ and some fixed Riemann-zeta zero 
$2z_0$ on the critical line. To this end, we employ a classical formula due to 
Potter [17], f.~(2.22), which approximates the Epstein zeta-function by a 
partial sum of its series, throughout the half-plane 
$\Re(s)>-{\text{1\over4}}$, $s\ne1$. In our notation, 
$$ \zeta_Q(s) =  F_1(Z,s) + F_2(Z,s)\,,  \eqno(4.3)  $$ 
$Z$ a positive real parameter, 
$$ \eq{F_1(Z,s) := & \sum_{\m\in\zzi\atop Q\m\le Z} Q\m^{-s} + 
s Z^{-s-1} \sum_{\m\in\Zi^2\atop Q\m\le Z} Q\m \ - \cr 
& - (1+s)Z^{-s} \sum_{\m\in\Zi^2\atop Q\m\le Z} 1 
+ {\pi\over\sqrt{D}}\,{s(s+1)\over(s-1)} Z^{1-s}\,,  \cr  } \eqno(4.4)  $$ 
$$ F_2(Z,s) := s(s+1) \Int_Z^\infty v^{-s-2} P_1(v) \d v\,, \eqno(4.5)  $$ 
where, for $v>0$, 
$$ P_1(v) := \Int_0^v P(w)\d w = {\sqrt{D}\over2\pi}\,v\,\sum_{\m\in\zzi} 
Q\m^{-1} J_2\(4\pi\sqrt{{v\over D} Q\m}\,\)\,,  $$ 
$J_2$ the usual Bessel function (see [17], Lemma 1). To estimate 
$\abs{F_2(Z,s)}$, we use that, for $x>0$, $\abs{J_2(x)}\le x^{-1/2}$, which is 
easily verified by f.~8.451 in [5]. This gives 
$$  \abs{P_1(v)} \le {D^{3/4} \over 4\pi^{3/2}}\,v^{3/4}\,
\sum_{\m\in\zzi} Q\m^{-5/4}\,. \eqno(4.6) $$ 
To bound this series, let $\disp 
\kappa_Q := \inf_{(u,v)\in\Ri_*^2} {Q(u,v)\over u^2+v^2}$, then a calculus 
exercise yields: If $\tau_\pm := {1\over b}(a-c\pm\sqrt{(a-c)^2+b^2})$ for 
$b\ne0$, then 
$$  \kappa_Q = \cases{\min\({Q(\tau_+,1)\over \tau_+^2+1},\ 
{Q(\tau_-,1)\over \tau_-^2+1}\) & if $b\ne0$, \cr \min(a,c) & if $b=0$.}  
\eqno(4.7)   $$ 
Hence 
$$  \sum_{\m\in\zzi} Q\m^{-5/4} \le \kappa_Q^{-5/4} \sum_{k=1}^\infty 
r(k) k^{-5/4} = 4 \kappa_Q^{-5/4} \zeta({\text{5\over4}}) 
L({\text{5\over4}})\,, $$ 
where $r(k)$ counts the number of ways to express $k$ as a sum of two squares, 
and $L(s)$ is the Dirichlet $L$-series\fn{9}{The evaluation of 
$L({\text{5\over4}})$ can be done by {\it Mathematica} [21], via the identity 
$L(s) = 2^{-s} \Phi(-1,s,\h)$, where $\Phi$ is the Lerch Phi-function: see 
[5], 
f.~9.550.} corresponding to the non-principal Dirichlet character 
\hbox{mod 4.} Combining this with (4.5) and (4.6), 
we obtain altogether, provided that\fn{10}{For better convergence, our 
strategy is to bound $|\zeta_Q(1-z_0)|$ away from 0, 
and then to appeal to the functional equation.} $\Re(s)={\text{3\over4}}$, 
$$  \abs{F_2(Z,s)} \le \abs{s(s+1)}\,{D^{3/4} \over \pi^{3/2}}\,
\kappa_Q^{-5/4} \zeta({\text{5\over4}}) L({\text{5\over4}})\,{1\over Z}\,. 
\eqno(4.8)  $$  
For a given form $Q$, one can therefore proceed as follows: Choose, e.g., 
$\z={\text{1\over4}}+i\b$ with $\b=7.06736\dots$ so that $\zeta(2\z)=0$, 
then by the functional equation (1.7), 
$$ \abs{\zeta_Q(\z)} = \({2\pi\over\sqrt{D}}\)^{-1/2}
{\abs{\Gamma(1-\z)}\over\abs{\Gamma(\z)}} \abs{\zeta_Q(1-\z)}\,.  $$ 
Combining this with (4.2), (3.8), and (4.3), we arrive at 
$$ \hbox{$ K_0 \ge 6\pi\,\abs{(\z-1)(2\z-1)\over(\z+2)^7}\,
\({2\pi\over\sqrt{D}}\)^{-1/2}
{\abs{\Gamma(1-\z)}\over\abs{\Gamma(\z)}}\,
\(\abs{F_1(Z,1-\z)}-\abs{F_2(Z,1-\z)}\)$}\,,  \eqno(4.9)  $$ 
where $Z$ remains a free parameter and $\abs{F_1(Z,1-\z)},\ \abs{F_2(Z,1-\z)}$ 
can be evaluated, resp., estimated by (4.4), (4.8). The only thing that could 
go wrong is that $\abs{\zeta_Q(1-\z)}$ is so small (or actually 0) that we 
cannot get a positive 
lower bound for the last bracket in (4.9). In this case, we can 
take one of the next Riemann-zeta zeros instead of $2\z$. \bsk 

{\bf Example.}\quad Let us consider the special (irrational) quadratic form 
$$ Q_0\m = m^2 + \sqrt{2}\, m n + \sqrt{3}\, n^2\,.  $$ 
Choosing $Z=1000$ and employing {\it Mathematica } [21] to evaluate (4.4), 
resp., (4.8), we obtain 
$\abs{F_1(1000,1-\z)}=0.422182\dots$, $\abs{F_2(1000,1-\z)}\le0.236529\dots$, 
hence $\abs{F_1(1000,1-\z)}-\abs{F_2(1000,1-\z)}\ge0.185653\dots$. Using this 
in (4.9), we finally arrive at 
$$ K_0 = \liminf_{Y\to\infty}\(Y^{-5/4} \Int_1^Y \abs{R(x)} \d x\) > 
4 \times 10^{-4}  $$ for this particular form $Q_0$. \msk 
Applying to the integral 
$\Int_{-1-i\infty}^{-1+i\infty} g(s) H(s) Y^s \d s$ in (3.7) similar 
arguments as we used to estimate $w(\eta)$, one can replace the 
$\liminf$ -bound by an inequality valid for all $Y>0$. For the form $Q_0$ we 
obtain in this way 
$$ Y^{-5/4} \Int_1^Y \abs{R(x)} \d x > 4 \times 10^{-4} - 3.62\,Y^{-5/4}\,, $$  
which is non-trivial for $Y>1500$.

\bsk 



\klein \parindent=0pt 

\vbox{\cen{\bf References}  \bsk 

[1] {\smc P.~Bleher,} 
On the distribution of the number of lattice points inside a 
family of convex ovals. Duke Math.~J. {\bf 67}, 461--481 (1992).  \ssk 

[2] {\smc J.B.~Conrey,} More than two fifth of the zeros of the Riemann 
zeta-function are on the critical line. 
J.~Reine Angew.~Math. {\bf 399}, 1--26 (1989). \ssk 
                                                        
[3] {\smc H.~Davenport \and H.~Heilbronn,} On the zeros of certain Dirichlet 
series I. J.~London Math.~Soc. {\bf 11}, 181--185 (1936). \ssk 

[4] {\smc H.~Davenport \and H.~Heilbronn,} On the zeros of certain Dirichlet 
series II. J.~London Math.~Soc. {\bf 11}, 307--312 (1936).} \ssk 

[5] {\smc I.S.~Gradshteyn \and I.M.~Ryzhik,} Table of integrals, series, and 
products. A.~Jeffrey editor. 5th ed., San Diego 1994.   \ssk 

[6] {\smc M.N.~Huxley}, {Exponential sums and lattice points II}. 
{Proc.~London Math.~Soc.} {\bf 66}, 279--301 (1993).  \ssk  
 
[7] {\smc M.N.~Huxley}, {Area, lattice points, and exponential sums.} 
LMS Monographs, New Ser. {\bf 13}, Oxford 1996.  \ssk 
 
[8] {\smc M.N.~Huxley \and W.G.~Nowak}, Primitive lattice points in convex 
planar domains. Acta Arithm. {\bf 76}, 271--283 (1996).  \ssk 

[9] {\smc A.~Ivi\'c,} The Riemann zeta-function. New York 1985. \ssk 

[10] {\smc E.~Kr\"atzel,} Lattice points. Berlin 1988. \ssk 

[11] {\smc E.~Kr\"atzel,} Analytische Funktionen in der Zahlentheorie. 
Stuttgart 2000. \ssk 

[12] {\smc N.~Levinson,} More than one third of the zeros of Riemann's 
zeta-function are on $\si=\h$. Adv.~Math. {\bf 13}, 383--436 (1974).  \ssk 

[13] {\smc W.~M\"uller}, Lattice points in convex planar domains: Power 
moments with an application to primitive lattice points. In: Proc.~Number 
Theory Conf.~held in Vienna 1996, 
W.G.~Nowak and J.~Schoi\ss engeier eds., Vienna 1996, pp.~189--199.  \ssk   

[14] {\smc W.G.~Nowak,} 
An $\Omega$-estimate for the lattice rest of a convex planar domain. 
Proc.~R.~Soc.~Edinb., Sect. A, {\bf 100}, 295-299 (1985).  \ssk 

[15] {\smc W.G.~Nowak,} On the mean lattice point discrepancy of a convex 
disc. Arch.~Math.~(Basel) {\bf 78}, 241--248 (2002).  \ssk 

[16] {\smc J.~Pintz,} On the distribution of square-free numbers. J.~London 
Math.~Soc.~(2) {\bf 28}, 401--405 (1983).  \ssk  

[17] {\smc H.S.A.~Potter,} Approximate equations for the Epstein 
zeta-function. Proc.~London Math.~Soc. (2) {\bf 36}, 501--515 (1934).  \ssk  

[18] {\smc A.~Selberg,} On the zeros of Riemann's zeta-function. Skr.~Norske 
Vid. Akad. Oslo, no.~{\bf10}, \hbox{59 p.} (1943). \ssk 

[19] {\smc E.C.~Titchmarsh,} The theory of the Riemann zeta-function. 
2nd ed., revised by D.R.~Heath-Brown. Oxford 1986.  \ssk 

[20] {\smc M.~Voronin,} On the zeros of zeta-functions of quadratic forms. 
Trudy Mat.~Inst.~Steklov {\bf 142}, 135--147 (1976).  \ssk  

[21] {\smc Wolfram} Research, Inc., Mathematica 4.1. Champaign 2001. \ssk 

[22] {\smc J.~Wu,} On the primitive circle problem. Monatsh.~f.~Math. 
{\bf135}, 69-81 (2002).  \ssk 

[23] {\smc W.~Zhai \and X.D.~Cao,} On the number of coprime integer pairs 
within a circle. Acta Arithm. {\bf 90}, 1-16 (1999).  \ssk 

[24] {\smc W.~Zhai,} On primitive lattice points in planar domains. Acta 
Arithm. {\bf109}, 1-26 (2003).  

\bsk\bsk


\parindent=1.5true cm

\vbox{Werner Georg Nowak

Institut f\"ur Mathematik u.~Ang.Stat. 

Universit\"at f\"ur Bodenkultur 

Peter Jordan-Stra\ss e 82 

A-1190 Wien, Austria \ssk 

E-mail: {nowak@mail.boku.ac.at} \ssk 

Web: http://www.boku.ac.at/math/nth.html}

\bye